\documentclass{amsart}
\usepackage{amssymb,graphicx,url,treesymb}
\usepackage[noadjust]{cite}

\usepackage[line,ps,color,poly,all,rotate]{xy}

\newcommand{\NN}{\mathbb N}

\newcommand{\ul}[1]{{\underline{#1}}}
\newcommand{\YY}{\mathbb{Y}}
\def\ls{\dashv }
\def\rs{\vdash }

\theoremstyle{plain}
\newtheorem{thm}{Theorem}[section]

\newtheorem{lem}[thm]{Lemma}
\newtheorem{prop}[thm]{Proposition}
\newtheorem{cor}[thm]{Corollary}
\theoremstyle{definition}
\newtheorem{defn}[thm]{Definition}

\begin{document}

\title{The arithmetic of trees}
\author{Adriano Bruno}
\address{Department of Mathematics and Statistics\\Lederle Graduate Research Tower\\ University of Massachusetts\\Amherst, MA 01003-9305}
\email{bruno@math.umass.edu}
\author{Dan Yasaki}
\address{Department of Mathematics and Statistics\\Lederle Graduate Research Tower\\ University of Massachusetts\\Amherst, MA 01003-9305}
\email{yasaki@math.umass.edu}
\date{}
\thanks{The original manuscript was prepared with the \AmS-\LaTeX\ macro
system and the \Xy-pic\ package.}
\begin{abstract}
The arithmetic of the natural numbers $\NN$ can be extended to arithmetic operations on planar binary trees.  This gives rise to a non-commutative arithmetic theory.  In this exposition, we describe this \emph{arithmetree}, first defined by Loday, and investigate prime trees. 
\end{abstract}
\maketitle

\bibliographystyle{../amsplain_initials}

\begin{section}{Introduction}\label{chapter:introduction}
J.-L. Loday recently published a paper \emph{Arithmetree} \cite{Lo}, in which he defines arithmetic operations on the set $\YY$ of groves of planar binary trees.  These operations extend the usual addition and multiplication on the natural numbers $\NN$ in the sense that there is an embedding $\NN \hookrightarrow \YY$, and the multiplication and addition he defines become the usual ones when restricted to $\NN$.  Loday's  reasons for introducing these notions have to do with intricate algebraic structures known as dendriform algebras  \cite{Lodend}. 

Since the arithmetic extends the usual operations on $\NN$, one can ask many of the same questions that arise in the natural numbers.  In this exposition, we examine notions of primality, specifically studying \emph{prime trees}.  We will see that all trees of prime degree must be prime, but many trees of composite degree are also prime.  One should not be misled by the idea that arithmetree is an extension of the usual arithmetic on $\NN$.  Indeed, away from the image of $\NN$ in $\YY$, the arithmetic operations $+$ and $\times$ are non-commutative.  Both operations are associative, but multiplication is only distributive on the left with respect to $+$.  In the end it is somewhat surprising that there is a very natural copy of $\NN$ inside $\YY$.  

The paper is organized as follows.  Sections~\ref{sec:defn}--\ref{sec:times} summarize without proofs the results that we need from \cite{Lo}.  Specifically, basic definitions are given in Section~\ref{sec:defn} to set notation.  The embedding  $\NN \hookrightarrow \YY$ is given in Section~\ref{sec:natural}, and Section~\ref{sec:basicoperations} discusses the basic operations on groves.  Sections~\ref{sec:add} and \ref{sec:times} define the arithmetic on $\YY$.  Finally, Section~\ref{sec:results} discusses some new results and Section~\ref{sec:final} gives a few final remarks. 

These results grew out of an REU project in the summer of 2007 at the University of Massachusetts at Amherst, and the authors thank them for their support.   The second author would like to thank Paul Gunnells for introducing him to this very interesting topic, as well as all the help with typesetting and computing.
\end{section}
\begin{section}{Background}\label{sec:defn}
In this section,  we give the basic definitions and set notation.  

\begin{defn}
A \emph{planar binary tree} is an oriented planar graph drawn in the plane with one root, $n+1$ leaves, and $n$ interior vertices, all of which are trivalent.\end{defn}  Henceforth, by \emph{tree}, we will mean a planar binary tree.  We consider trees to be the same if they can be moved in the plane to each other.  Thus we can always represent a tree by drawing a root and then having it
``grow'' upward.  The \emph{degree} is the number of
internal vertices.  See Figure \ref{fig:tree} for an example of a tree of degree four.
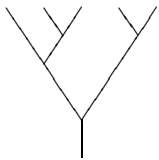
\begin{figure}[htb]
\[
\begin{xy}
(0,5)*{}="a1";(5,5)*{}="a2";(10,5)*{}="a3";(15,5)*{}="a4";(20,5)*{}="a5";
(7.5,1.25)*{}="b1";(17.5,1.25)*{}="b2";
(5,-2.5)*{}="c";
(10,-10)*{}="d";
(10,-15)*{}="r";
"a1";"d"**\dir{-};
"a2";"b1"**\dir{-};
"a3";"c"**\dir{-};
"a4";"b2"**\dir{-};
"a5";"d"**\dir{-};
"d";"r"**\dir{-};
\end{xy}
\]
\caption{\label{fig:tree}A tree of degree four with five leaves.}
\end{figure}

Let $Y_{n}$ be the set of trees of degree $n$.  For example,
\begin{gather*}
Y_{0} = \{\one \},\quad
Y_{1} = \{\oneone \},\quad
Y_{2} = \{\twoone, \onetwo \},\quad\text {and}\quad Y_{3} = \{\threeone, \onethree, \twotwo, \onetwoonea, \onetwooneb \}.
\end{gather*}
One can show that the cardinality of $Y_{n}$ is given by the $n^{th}$ \emph{Catalan} number,  
\[c_n= \frac{1}{n+1} \binom{2n}{n}=\frac{(2n)!}{(n+1)!n!}.\] 
The Catalan numbers arise in a variety of combinatorial problems \cite{Stan}.\footnote{He currently gives 161 combinatorial interpretations of $c_n$.} 
\begin{defn}
A nonempty subset of $Y_{n}$ is called a \emph{grove}.  The set of all
groves of degree $n$ is denoted by $\YY_{n}$. \end{defn} For example,
\[
\YY_{0} = \{\one \}, \quad
\YY_{1} =\{\oneone \}, \quad \text{and}\quad  \YY_{2} = \{\twoone, \onetwo, \twoone\cup \onetwo \}.
\]
Notice that we are omitting the braces around the sets in $\YY_{n}$ and use instead $\cup$ to denote the subsets.  For example we write $\twoone\cup \onetwo$ as opposed to $\{\twoone,\onetwo\}$ to denote the grove in $\YY_2$ consisting of both trees of degree $2$.  Let $\YY=\bigcup_{n \in \NN} \YY_n$ denote the set of all groves.  By definition groves consist of trees of the same degree; hence we get a well-defined notion of degree \begin{equation}\deg:\YY\to \NN.\end{equation}

The Catalan numbers $c_n$ grow rapidly.  Since $\YY_n$ is the set of subsets of $Y_n$, we see that the cardinality $\#\YY_n = 2^{c_n}-1$ grows extremely fast, necessitating the use of computers even for computations on trees of fairly small degree.
\begin{table}[htb]
\caption{\label{tab:size}Number of trees and groves of degree $n \leq 7$.}
\begin{tabular}{|r|r|r|}
\hline
$n$&$\#Y_n$&$\#\YY_n$\\\hline
1& 1& 1\\
2& 2& 3\\
3& 5& 31\\
4& 14& 16383\\
5& 42& 4398046511103\\
6& 132& 5444517870735015415413993718908291383295\\
7& 429& $\sim 1.386 \times 10^{129}$\\\hline
\end{tabular}
\end{table}

\end{section}
\begin{section}{The natural numbers}\label{sec:natural}
In this section we give an embedding of $\NN$ into $\YY$.  There is a distinguished grove for each degree given by set of all trees of degree $n$.
\begin{defn}
The \emph{total grove of degree $n$} is defined by $\ul{n}=\bigcup_{x\in Y_{n}} x$.
\end{defn}
For example
\[\ul{0}=\one, \quad \ul{1}=\oneone, \quad \ul{2}=\onetwo \cup \twoone,\quad \text{and}  \quad \ul{3}=\threeone \cup \onetwoonea \cup \onetwooneb \cup \onethree \cup \twotwo \]
This gives an embedding $\NN \hookrightarrow \YY$.  It is clear that the degree map is a one-sided inverse in the sense that $\deg(\ul{n})=n$ for all $n \in \NN$.  We will see in Section~\ref{sec:properties} that under this embedding, \emph{arithmetree} can be viewed as an extension of arithmetic on $\NN$.
\end{section}
\begin{section}{Basic operations}\label{sec:basicoperations}
In this section we define a few operations that will be used to define the arithmetic on $\YY$.  
\subsection{Grafting}
\begin{defn}
We say that a tree $z$ is obtain as the \emph{graft} of $x$ and
$y$ (notation: $z=x\vee y$) if $z$ is gotten by attaching the root of $x$ to the
left leaf and the root of $y$ to the right leaf of $\oneone$.  
\end{defn} For example, $\twoone = \oneone \vee \one$ and $\twotwo = \oneone \vee \oneone$. 
It is clear that every tree $x$ of degree greater than $1$ can be obtained as the graft of trees $x^l$ and $x^r$ of degree less than $n$.  Specifically, we have that $x=x^l \vee x^r$.  We refer to these subtrees as the \emph{left} and \emph{right parts} of $x$.

Given a tree $x$ of degree $n$, then one can create a tree of degree $n+1$ that carries much of the structure of $x$ by grafting on $\ul{0}=\one$.  Indeed, there are two such trees, $x \vee \ul{0}$ and $\ul{0} \vee x$.  We will say that such trees are \emph{inherited}.
\begin{defn}
A tree $x$ is said to be \emph{left-inherited} if $x^r=\ul{0}$ and \emph{right-inherited} if $x^l=\ul{0}$.  A grove is \emph{left-inherited} (resp. \emph{right-inherited}) if each of its member trees is \emph{left-inherited} (resp. \emph{right-inherited}).   
\end{defn}

We single out two special sequences of trees $L_n$ and $R_n$.  
\begin{defn}
Let $L_1=R_1=\ul{1}$.  For $n>1$, set $L_n=L_{n-1}\vee \ul{0}$ and $R_n=\ul{0} \vee R_{n-1}$.  We will call such trees \emph{primitive}.
\end{defn}
Notice that $L_n$ is the left-inherited tree such that $L_n^l=L_{n-1}$.  Similarly,  $R_n$ is the right-inherited tree such that $R_n^r=R_{n-1}$.  
\subsection{Over and under}
\begin{defn}
For $x\in Y_p$ and $y\in Y_q$ the tree $x/y$
(read $x$ {\it over} $y$) in $Y_{p+q}$ is obtained by identifying the
root of $x$ with the leftmost leaf of $y$.  Similarly, the tree
$x\backslash y$ (read $x$ {\it under} $y$) in $Y_{p+q}$ is obtained by
identifying the rightmost leaf of $x$ with the root of $y$.
\end{defn}
For example, $\onetwo / \oneone = \onetwoonea$ and $\twoone \backslash
\oneone = \twotwo$.
\subsection{Involution}
The symmetry around the axis passing through the root defines an involution $\sigma$ on $Y$.  For example, $\sigma(\twotwo)=\twotwo$ and $\sigma(\twoone)=\onetwo$.  The involution can be extended to an involution on $\YY$, by letting $\sigma$ act on each tree in the grove.  On can easily check that for trees $x,y$,
\begin{enumerate}
\item $\sigma(x \vee y)=\sigma(y)\vee \sigma(x)$,
\item $\sigma(x/y)=\sigma(y)\backslash \sigma(x)$, and 
\item $\sigma(x \backslash y)=\sigma(y)/\sigma(x)$.
\end{enumerate}
We will see that this involution also respects the arithmetic of groves.
\end{section}
\begin{section}{Addition}\label{sec:add}
Before we define addition, we first put a partial ordering on $Y_n$.  
\subsection{Partial ordering}
We say that the inequality $x<y$ holds if $y$ is obtained from
$x$ by moving edges of $x$ from left to right over a vertex.  This induces a partial ordering on $Y_n$ by imposing
\begin{enumerate}
\item $(x \vee y)\vee z \leq x \vee (y\vee z)$
\item If $x < y$ then $x \vee z < y \vee z$  and $z \vee x < z \vee y$ for all $z \in Y_n$. 
\end{enumerate}
For example,  $\threeone < \onetwoonea < \onetwooneb < \onethree$.  Note that the primitive trees are extremal elements with respect to this ordering.
\subsection{Sum}
\begin{defn}
The {\it sum} of two trees $x$ and
$y$ is the following disjoint union of trees
\[ x + y := \bigcup_{x/y \le z \le x\backslash  y} z\ . \]
\end{defn}
All the elements in the sum have the same  degree which happens to be $\deg(x) + \deg(y)$.  Thus we can extend the definition of addition to groves by distributing.  Namely, for groves  $x=\bigcup_i x_i$ and $y=\bigcup_j y_i$, 
\begin{equation}\label{eq:addunion}
x + y := \bigcup_{ij}\, (x_i+y_j).\end{equation}
\begin{prop}[Recursive property of addition]\label{prop:recursiveadd}
Let  $x=x^{l}\vee x^{r}$ and $y=y^{l}\vee y^{r}$ be non-zero trees.  Then 
\[
x + y = x^l\vee (x^r + y)\, \cup \, (x+y^l)\vee y^r.
\]
\end{prop}

Note that the recursive property of addition says that the sum of two trees $x$ and $y$ is naturally a union of two sets, which we call the \emph{left} and \emph{right sum} of $x$ and $y$:
\begin{equation}\label{eq:lsrs}
x\ls y = x^l\vee (x^r + y)\quad \text{and} \quad x \rs y =(x+y^l)\vee y^r.\footnote{We set $x \rs \ul{0} = \ul{0} \ls y = \ul{0}$.}
\end{equation}
Note that $x+y = x\ls y \cup x \rs y$.  You can think about this as
splitting the plus sign $+$ into two signs $\ls$ and $\rs$.  From \eqref{eq:addunion} and the definition, we see that the definition for left sum and right sum can also be extended to groves by distributing.

With the definition of inherited trees/groves and \eqref{eq:lsrs}, one can easily check that left (respectively right) inheritance is passed along via right (respectively left) sums.  More precisely,
\begin{lem}\label{lem:inherit}
Let $y$ be a left-inherited tree.  Then $x \rs y$ is left-inherited.  Similarly, if $x$ is right-inherited, then $x \ls y$ is right-inherited.
\end{lem}
\subsection{Universal expression}
It turns out that every tree can expressed as a combination of left and right sums of $\oneone$.  This expression is unique modulo the failure of left and right sum to be associative.  More precisely, 
\begin{prop}
Every tree $x$ of degree $n$ can be written in as an iterated Left
and Right sum of $n$ copies of $\oneone$.  This is called the
\emph{universal expression} of $x$, and we denote it by $w_{x}
(\oneone)$.  
This expression is unique modulo 
\begin{enumerate}
\item $(x \ls y) \ls z = x \ls (y+z)$,
\item $(x \rs y) \ls z = x \rs (y \ls z)$, and 
\item $(x+y) \rs z = x \rs (y \rs z)$.
\end{enumerate}
\end{prop}
For example
\[\twoone = \oneone\rs\oneone \quad \twotwo = \oneone \rs \oneone \ls \oneone.
\]  
Loday gives a algorithm for computing the universal expression of a tree $x$.

\begin{prop}[Recursive property for universal expression]\label{prop:universal}
Let $x$ be a tree of degree greater than $1$.  The algorithm for determining $w_{x} (\oneone)$ is given through the
recursive relation
\[
w_{x} (\oneone) = w_{x^{l}} (\oneone) \rs \oneone \ls w_{x^{r}} (\oneone).
\]
\end{prop}
\end{section}

\begin{section}{Multiplication}\label{sec:times}
Essentially, we define the multiplication to distribute on the left over the universal expression.  
\begin{defn}
The product $x\times y$ is
defined by 
\[
x\times y = w_{x} (y).
\]
\end{defn}
This means to compute the product $x\times y$, first compute the
universal expression for $x$, then replace each occurrence of $\oneone$
by the tree $y$, then compute the resulting Left and Right sums.  For example, one can easily check that $\twoone = \oneone\rs\oneone$.  This means for any tree $y$, $\twoone \times y = y \rs y$.  In particular
$\twoone\times \onetwo = \onetwo \rs \onetwo$ is the tree shown in Figure~\ref{fig:tree}.

Note that the definition of $x \times y$ as stated still makes sense if $y$ is a grove.  We can further extend the definition of multiplication to the case when $x$ is a grove by declaring multiplication to be distributive on the left over disjoint unions: \[(x \cup x') \times y = x \times y \cup x' \times y =w_x(y) \cup w_{x'}(y).\]
\end{section}

\begin{section}{Properties}\label{sec:properties}
We list a few properties of \emph{arithmetree}.
\begin{itemize}
\item The addition $+: \YY \times \YY \to \YY$ is associative, but not commutative. 
\item The multiplication $\times: \YY \times \YY \to \YY$ is associative, but not commutative.  It is distributive on the left with respect to $+$, but it is not right distributive. 
\item There is an injective map $\NN \hookrightarrow \YY$, $n \mapsto \ul{n}$ (defined in Section~\ref{sec:natural}) that respects the arithmetic.  Namely, 
\[\ul{m+n}=\ul{m}+\ul{n}\quad \text{and}\quad \ul{mn}=\ul{m}\times\ul{n} \quad \text{for all $m,n \in \NN$.}\]
\item Degree gives a surjective map $\deg:\YY \to \NN$ that respects the arithmetic and is a one-sided inverse to the injection above .  For every $x,y \in \YY$, 
\[
\deg(x+y)=\deg(x)+\deg(y) \quad \text{and}
\quad \deg(x\times y)=\deg(x)\deg(y).\]
\item $\deg(\ul{n})=n$ for all $n \in \NN$.
\item The neutral element for $+$ is $\ul{0}=\one$.
\item The neutral element for $\times$ is $\ul{1}=\oneone$.
\item The involution $\sigma$ satisfies
\[\sigma(x+y)=\sigma(y)+\sigma(x) \quad \text{and} \quad \sigma(x\times y)=\sigma(x) \times \sigma(y).\]
\end{itemize}
\end{section}
\begin{section}{Results}\label{sec:results}
The recursive properties of addition and multiplication allowed us to implement arithmetree on a computer using gp/PARI \cite{gp}.  The computational experimentation was done using Loday's naming convention for trees \cite{Lo}.
\subsection{Counting trees}
Since each grove $x \in \YY$ is just a subset of trees, there is another measure of the ``size'' of $x$ other than degree. 
\begin{defn}
Let $x \in \YY$ be a grove.  The \emph{count} of $x$, denoted $C(x)$ is defined as the cardinality of $x$.
\end{defn}
It turns out that count function gives a coarse measure of how complicated a grove $x$ is in terms of arithmetree.  Namely, if $x$ is the sum (resp. product) of other groves, then the count of $x$ is at least as large as the count of any of the summands (resp. factors).
\begin{lem}\label{lem:lsrs}
Let $x ,y \in \YY$ be two non-zero groves.  Then 
\begin{enumerate}
\item $C(x \ls y) \geq C(x)C(y)$, with equality if and only if $x$ is a left-inherited grove. \label{it:ls}
\item $C(x \rs y) \geq C(x)C(y)$, with equality if and only if $y$ is a right-inherited grove.\label{it:rs}
\end{enumerate}
\end{lem}
\begin{proof}
We first consider \eqref{it:ls}.  Since $\ls$ is distributive over unions, it suffices to prove the case when $x$ and $y$ are trees.  Namely, we must show that for all non-zero trees $x$ and $y$, $C(x \ls y) \geq 1$, with equality if and only if $x$ is a left-inherited tree. It is immediate that $C(x \ls y) \geq 1$; it remains to show that equality is only attained when $x$ is left-inherited.  From the definition of left sum, $x \ls y =x^l \vee (x^r + y)$.  If $x$ is not left-inherited, then $x^r \neq \ul{0}$ and  
\begin{align*}
C(x \ls y) &= C(x^l \vee (x^r + y))\\
&=C(x^r + y)\\
&=C(x^r\ls y \cup x^r \rs y)\\
&=C(x^r\ls y)+ C(x^r \rs y)\\
&>1.
\end{align*}
On the other hand, if $x$ is left-inherited, then $x^r = \ul{0}$ and  
\[C(x \ls y) = C(x^l \vee (x^r + y))=C(x^l\vee y)=1.\]

Item \eqref{it:rs} follows similarly.
\end{proof}
\begin{prop}\label{prop:sp}
Let $x ,y \in \YY$ be two non-zero groves.  Then 
\begin{enumerate}
\item $C(x + y) \geq 2C(x)C(y)$, with equality if and only if $x$ is a left-inherited and $y$ is right-inherited.\label{it:s}
\item $C(x \times y) \geq C(x)C(y)^{\deg(x)}$. \label{it:p}
\end{enumerate}
\end{prop}
\begin{proof}
Since $x+y=x \ls y \cup x\rs y$, \eqref{it:s} follows immediately from Lemma~\ref{lem:lsrs}.  For \eqref{it:p}, we note that multiplication is left distributive over unions, and so it suffices to prove the case when $x$ is a tree.  Namely we must show that for a tree $x$ and a grove $y$, $C(x \times y) \geq C(y)^{\deg(x)}$.

Let $w_x$ be the universal expression of the tree $x$.  Then $x \times y=w_x(y)$ is some combination of left and right sums of $y$.  By distributivity of left and right sum over unions and repeated usage of Lemma~\ref{lem:lsrs}, the result follows.
\end{proof}
\subsection{Primes}
\begin{defn}
A grove $x$ is said to be \emph{prime} if $x$ is not the product of two groves different from $\ul{1}$. 
\end{defn}
Since $\deg(x \times y)=\deg(x)\deg(y)$ for all groves $x,y$, it is immediate that any grove of prime degree is prime.  However, there are also prime groves of composite degree.  For example, by taking all possible products of elements of $\YY_2$, one can check by hand that the primitive tree $L_4$ is a prime grove of degree $4$.  

We turn our focus to prime trees, which are prime groves with count equal to $1$.  It turn out that composite trees have a nice description in terms of inherited trees.  Namely, a composite tree must have an inherited tree as a right factor and a primitive tree as a left factor.
\begin{thm}\label{thm:factor}
Let $z$ be a composite tree of degree $n$.  Then there exists a proper divisor $d\neq 1$ of $n$ and a tree $T \in Y_{d-1}$ such that 
\[z=L_{n/d} \times (\ul{0} \vee T) \quad \text{or} \quad z=R_{n/d} \times (T \vee \ul{0})\]
\end{thm}    
\begin{proof}
Let $z=x\times y$ be a composite tree of degree $n$.  By Proposition~\ref{prop:sp}, $x$ and $y$ must also be trees.  Since $n=\deg(z)=\deg(x)\deg(y)$, it follows that there exists a proper divisor $d\neq 1$ of $n$ such that $\deg(y)=d$ and $\deg(x)=n/d$.

We proceed by induction on the degree of $x$.  Suppose $x$ is a tree of degree $2$.  Then $x=\oneone\ls \oneone$ or $x=\oneone\rs \oneone$.  If $x=\oneone\rs \oneone$, then $x=L_2$ is primitive and  
\[1=C(x\times y) =C(y \rs y).\]
From Proposition~\ref{prop:sp}, it follows that $y$ is $y$ is right-inherited.  Similarly, if $x=\one \ls \one$, then $x=R_2$ and $y$ is left-inherited.

Now suppose $x$ is a tree of degree $k$ such that $x \times y$ is a tree of degree $n$.  From Proposition~\ref{prop:universal} and the definition of multiplication, it follows that 
\begin{align*}
x \times y &=w_x(y)\\
&=w_{x^l}(y) \rs y \ls w_{x^r}(y)\\
&=(x^l \times y) \rs y \ls (x^r \times y).
\end{align*}

Suppose $x^r \neq \ul{0}$.  Then  $x^r \times y \neq \ul{0}$ and $C(y \ls (x^r \times y))=1$.  Then by Proposition~\ref{prop:sp}, $y$ is left-inherited.  Let $T=y \ls (x^r \times y)$.  By Lemma~\ref{lem:inherit}, $T$ is also left-inherited.  Since $C((x^l \times y) \rs T)=1$ and $T \neq \ul{0}$, we must have that either $T$ is also right-inherited, or $(x^l \times y)=\ul{0}$.  The only tree that is both left and right-inherited is the tree $\ul{1}=\one$.  It follows that  $(x^l \times y)=\ul{0}$, and hence $x^l=\ul{0}$.  By the inductive hypothesis, $x^r$ is a right-primitive tree, and hence $x=R_k$.   

Now suppose $x^r =\ul{0}$.  Then $x^l \neq \ul{0}$, and an analogous argument shows that $y$ is left-inherited and $x=L_k$.   
\end{proof}

From this theorem, we get a nice picture of composite trees as the trees with a particular shape.  More precisely, one computes that the products $L_k \times (\ul{0} \vee T)$ and $R_k \times (T \vee \ul{0})$ have the forms given in Figure~\ref{fig:composite}.  It follows that the primitive trees ($L_k$ and $R_k$) and the inherited trees ($\ul{0}\vee T$ and $T \vee \ul{0}$) are prime.  More precisely, 
\begin{prop}\label{prop:unique}
A non-zero tree is either $\oneone$, prime, or the product of exactly two prime trees.  Furthermore, the factors are exactly the ones given in Theorem~\ref{thm:factor}, and can be read off from the shape of the tree. 
\end{prop}

\begin{figure}
\[
\begin{array}{c@{\hspace{0.4in}} c}
\begin{xy}
(0,5)*{}="a1";(5,5)*{T}="a2";(10,5)*{\cdots}="a3";(15,5)*{T}="a4";(20,5)*{T}="a5";(25,5)*{T}="a6";
(2.5,1.25)*{}="b1";
(7.5,-6.25)*{}="c";
(10,-10)*{}="d";
(12.5,-13.75)*{}="e";
(12.5,-17.5)*{}="r";
"a2";"b1"**\dir{-};
"a5";"d"**\dir{-};
"a4";"c"**\dir{-};
"e";"r"**\dir{-};
"e";"a1"**\dir{-};
"e";"a6"**\dir{-};
\end{xy} & \begin{xy}
(0,5)*{T}="a1";(5,5)*{T}="a2";(10,5)*{T}="a3";(15,5)*{\cdots}="a4";(20,5)*{T}="a5";(25,5)*{}="a6";
(22.5,1.25)*{}="b1";
(17.5,-6.25)*{}="c";
(15,-10)*{}="d";
(12.5,-13.75)*{}="e";
(12.5,-17.5)*{}="r";
"a5";"b1"**\dir{-};
"a2";"d"**\dir{-};
"a3";"c"**\dir{-};
"e";"r"**\dir{-};
"e";"a1"**\dir{-};
"e";"a6"**\dir{-};
\end{xy}

\end{array}
\]
\caption{\label{fig:composite} These are the possible shapes of composite trees.  The figure on the left is $L_k \times (\ul{0} \vee T)$ and the one on the right is $R_k \times (T\vee \ul{0})$.}
\end{figure}
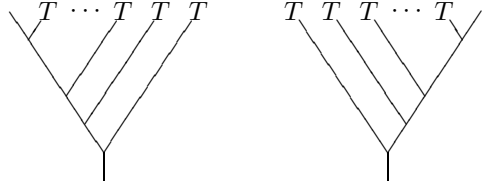
  
As a consequence of Proposition~\ref{prop:unique}, one has the following combinatorial formula.
\begin{cor}
Let $a_n$ denote the number of composite trees of degree $n$.  Then 
\[\frac{a_n}{2} = -c_1-c_n+\sum_{d \mid n}c_d,\quad \text{where $c_d$ is the $d^{th}$ Catalan number}.\]
\end{cor}
\end{section}
\begin{section}{Final remarks}\label{sec:final}
\subsection{Unique factorization}
In \cite{Lo}, Loday conjectures that arithmetree possesses unique factorization.  Namely, when a grove $x$ is written as a product of prime groves, the ordered sequence of factors is unique.  Very narrowly interpreted, this statement is false.  For example since multiplication in $\NN$ is commutative and multiplication in $\YY$ extends arithmetic on $\NN$, we see that for $n\in \NN$, if $n=p_1p_2 \cdots p_k$, then $\ul{n}=\ul{p_{\sigma(1)}}\times\ul{p_{\sigma(2)}}\times \cdots \times \ul{p_{\sigma(k)}}$ for any permutation $\sigma$.  However, away from the image of $\NN$ in $\YY$, it appears that this narrow interpretation is true.  Specifically, computer experimentation on groves of degree up to $12$ yielded a unique ordered sequence of prime factors for each grove outside of the image of $\NN$ in $\YY$.  

If we interpret the image of $\NN$ in $\YY$ in terms of the count function, we see that it is precisely the set of groves with maximal count;
\[\YY^{\max}=\bigcup_{n \in \NN}\{x \in \YY_n\;|\;C(x)=c_n\}.\]  
This subset $\YY^{\max}$ possesses unique factorization up to permutation of the factors.  On the other extreme, the trees are precisely the set of groves with minimal count;
\[\YY^{\min}=\bigcup_{n \in \NN}\{x \in \YY_n\;|\;C(x)=1\}.\]  
It follows from Proposition~\ref{prop:unique} that $\YY^{\min}$ possesses unique factorization in the narrow sense. The question of unique factorization for all of $\YY$ is open. 
\subsection{Additively irreducible}
From Proposition~\ref{prop:sp} we see that not every grove can be written as a sum of groves.  In fact it is easy to see that every tree is \emph{additively irreducible} in the sense that it cannot be written as the sum of two groves.  It would be interesting to study  additively irreducible groves.  In an analogue to the question of unique factorization, one could ask if arithmetree possesses \emph{unique partitioning}.  Namely, when a grove is written as a sum of additively irreducible elements, is the ordered sequence of summands unique?  
\end{section}
\bibliography{../references}    
\end{document}